\documentclass[journal]{IEEEtran}
\usepackage{amssymb}
\usepackage{epsfig}
\usepackage{graphicx}
\usepackage{colortbl}
\usepackage{float}
\usepackage{multirow}

\newtheorem{theorem}{Theorem}

\newtheorem{lemma}{Lemma}

\newtheorem{definition}{Definition}

\usepackage{algorithm}
\usepackage{algorithmic}
\ifCLASSINFOpdf
\else
\fi
\begin{document}
%
\title{Q-Learning for Linear Quadratic Optimal Control with Terminal State Constraint }
%
%
%

\author{Juanjuan~Xu,~Jingmei~Liu,~Zhaorong~Zhang,~Wei~Wang
\thanks{*This work was supported by the National Natural Science Foundation of
China under Grants 61821004, 62250056, and the Natural Science Foundation of Shandong Province (ZR2021ZD14, ZR2021JQ24),
and High-level Talent Team Project of Qingdao West Coast New Area (RCTDJC-2019-05), Key Research and Development Program of Shandong Province
(2020CXGC01208), and Science and Technology Project of Qingdao West
Coast New Area (2019-32, 2020-20, 2020-1-4). W. Wang is the corresponding author.}
\thanks{J. Xu is with School of Control Science and Engineering, Shandong University, Jinan, Shandong, P.R. China 250061.
        {\tt\small juanjuanxu@sdu.edu.cn}}
\thanks{J. Liu is with School of Control Science and Engineering, Shandong University, Jinan, Shandong, P.R. China 250061.
        {\tt\small jingmei0729@qq.com}}
\thanks{Z. Zhang is with School of Computer Science and Technology,Shandong University, Qingdao, Shandong, China 266237.
        {\tt\small zhaorong.zhang@uon.edu.au}}
\thanks{W. Wang is with School of Control
Science and Engineering, Shandong University, Jinan, Shandong, P.R.China 250061.
        {\tt\small w.wang@sdu.edu.cn}}%
}

\maketitle

\begin{abstract}
This paper is concerned with the linear quadratic optimal control of discrete-time time-varying system with terminal state constraint.
The main contribution is to propose a Q-learning algorithm for the optimal controller when the time-varying system matrices and input matrices are both unknown. Different from the existing Q-learning algorithms in the literature which are mainly for the unconstrained optimal control problem, the novelty of the proposed algorithm is available to deal with the case with terminal state constraints. A numerical example is illustrated to verify the effectiveness of the proposed algorithm.

\end{abstract}


%
\IEEEpeerreviewmaketitle

\section{Introduction}

Linear quadratic (LQ) optimal control has wide applications in lots of
fields including aerospace engineering \cite{background1},
economics \cite{background2} and so on, the research of which began in the 1950s \cite{background3}. Since then, plenty of essential progress has been made
\cite{progress1, progress2}. In \cite{Kalman}, the optimal solution was given in the feedback of the state in terms of Riccati equation.
\cite{Ignaciuk} provided a closed-form expression for the LQ optimal control of periodic-review perishable inventory systems.
\cite{Jiongmin} proved the existence of an equilibrium control for the deterministic LQ time-inconsistent optimal control problem.
\cite{Irregular} obtained the necessary and sufficient condition for the solvability of the LQ optimal control with irregular performance. These results are mainly derived for the optimal control problem without constraints.

Considering the practical demands in applications such as spacecraft launch and  mean-variance portfolio selection, the state/input must satisfy specific constraints, e.g., the terminal state is fixed \cite{demand1}, the control input is bounded \cite{demand2} etc. This stimulates the study of the LQ optimal control with constraints.
\cite{dp2020} derived the solution of a constrained LQ regulator problem via the set of optimal
active sets by using dynamic programming.
A constraint-compliant feedforward-feedback control update was derived in \cite{projection2017} for the equality constrained iterative LQ optimal control.
A Riccati-based approach was proposed in \cite{Riccati2010} to solve the LQ optimal control subject to linear
equality path constraints including mixed state-control and
state-only constraints.
In \cite{unified2006}, a closed-form expression parametrizing all the solutions of the
Hamiltonian system over a finite time interval was given by a strongly unmixed solution of a continuous algebraic Riccati equation and of the solution of an algebraic Lyapunov equation. LQ optimal control with fixed terminal states
and integral quadratic constraints was studied in \cite{sunjingrui}.
LQ tracking control with fixed terminal states was given in
recursive forms in \cite{Park2008}. As a brief summary, we note that most results for the constrained LQ optimal control is model-based, that is, the design of the optimal controller depends strongly on the exact knowledge of system matrices and input matrices.
However, a large number of engineering practices are faced with the uncertain environment and external
disturbances which makes it difficult to acquire the accurate models \cite{ZhangZhaorong}. Therefore, it is necessary to study the constrained
LQ optimal control with unknown system matrices and input matrices.

As is well known, Q-learning algorithm is a kind of reinforcement learning methods, which is very effective for dealing with the optimal control problem with unknown model. \cite{FinitehorizonLQ} presented a finite-horizon Q-learning algorithm which reveals state-feedback and output-feedback policies of standard LQ optimal control problem. \cite{ZhangZhaorong} proposed a Q-learning algorithm for the feedback Nash strategy of the finite-horizon nonzero-sum difference game.  \cite{LinZongli} provided the optimal
control policy for infinite-horizon LQ zero-sum game by utilizing a novel Q-Learning output feedback algorithm. In \cite{TaoFeng}, Q-learning algorithm was employed in solving consensusability problem of discrete-time multi-agent systems. It is worth noting that although the Q-learning algorithm has been widely studied in various problems and many important progress has been made,
the Q-learning for the LQ optimal control with terminal state constraints remains to be solved. The main difficulty comes from the learning of a specific parameter in view of the fact that the optimal control is in the feedback of the state and the specific parameter.

In this paper, we study the LQ optimal control of discrete-time time-varying system with terminal state constraint.
A Q-learning algorithm is proposed for the optimal controller when the time-varying system matrices and input matrices are both unknown. In particular, the Q-learning algorithm is consisting of two steps. The first step is to learn the solution of the Riccati equation, which thus gives the matrices gains of a specific matrix equation. The second step is to calculate the specific parameter from the data-driven matrix equation. Accordingly, the optimal controller is derived by running the Q-learning algorithm. Finally, we illustrate a numerical example to verify the effectiveness of the proposed algorithm.

The remainder of the paper is organized as follows. Section II presents the studied problem.
The model-based solution is solution in Section III. Section IV presents the detailed derivation of the Q-learning algorithm.
Numerical examples are given in Section V. Some concluding remarks are given in the last section.

The following notations will be used throughout this paper: $R^n$
denotes the set of $n$-dimensional vectors; $x'$ denotes the
transpose of $x$; a symmetric matrix $M>0\ (\geq 0)$ means that $M$ is
strictly positive-definite (positive semi-definite). $I$ denotes the identity matrix with appropriate dimension.


\section{Problem Formulation}

In this paper, we consider the following discrete-time time-varying system:
\begin{eqnarray}
x(k+1)&=&A(k)x(k)+B(k)u(k),\label{t1}
\end{eqnarray}
where $x\in R^n$ is the state, $u\in R^m$ is the control input.
$A(k)$ and $B(k)$ are time-varying matrices with compatible dimensions.
$k$ is the discrete-time variable taking values in $\{0, 1, \ldots, N\}$.
The initial value is prescribed as $x(0)=x$.

The cost function is given by
\begin{eqnarray}
J&=&\sum_{k=0}^N[x'(k)Qx(k)+u'(k)Ru(k)]\nonumber\\
&&+x'(N+1)Hx(N+1),\label{t2}
\end{eqnarray}
where $Q$ and $H$ are positive semi-definite matrices with compatible dimensions, and $R$ is a positive-definite matrix with compatible dimension.

The studied problem is presented as follows.

\textbf{Problem.} The aim of the paper is to design a Q-learning algorithm for the optimal controller $u(k)$ which minimizes the cost function (\ref{t2}) subject to system (\ref{t1}) and ensures the terminal state satisfies $x(N+1)=\xi$ with unknown model matrices $A(k), B(k), k=0, 1, \ldots, N$ where $\xi$ is a given constant vector.

\section{Model-based solution}

In this section, we first present the model-based optimal controller for the Problem. To this end, we introduce the Riccati difference equation:
\begin{eqnarray}
P(k)&=&A'(k)P(k+1)A(k)+Q-A'(k)P(k+1)B(k)\nonumber\\
&&\times [R+B'(k)P(k+1)(k)]^{-1}\nonumber\\
&&\times B'(k)P(k+1)A(k),\label{t3}
\end{eqnarray}
with terminal condition $P(N+1)=H$ and make the following denotations:
\begin{eqnarray}
\Gamma(k)&=&R+B'(k)P(k+1)B(k),\nonumber\\
K(k)&=&-\Gamma^{-1}(k)B'(k)P(k+1)A(k),\nonumber\\
A_c(k)&=&A(k)+B(k)K(k),\nonumber\\
\Phi(k,N)&=&\left\{
              \begin{array}{ll}
                A_c(N)\cdots A_c(k), & \hbox{$ k\leq N$;} \\
                I, & \hbox{$k>N$.}
              \end{array}
            \right.,\nonumber\\
K_1(k)&=&-\Gamma^{-1}(k)B'(k)\Phi'(k+1,N),\nonumber\\
\bar{B}(k)&=&B(k)\Gamma^{-1}(k) B'(k),\nonumber\\
G(s)&=&\sum_{j=s}^N\Phi(j+1,N)\bar{B}(j)\Phi'(j+1,N).\nonumber
\end{eqnarray}

\begin{lemma}\label{lem1}
Assume that the following equation
\begin{eqnarray}
\Phi(0,N)x-G(0)\lambda=\xi,\label{t4}
\end{eqnarray}
has a solution $\lambda^*$, then the optimal controller for the Problem is given by
\begin{eqnarray}
u(k)&=&K(k)x(k)+K_1(k)\lambda^*.\label{t5}
\end{eqnarray}
\end{lemma}
\emph{Proof.} Firstly, considering the terminal state constraint $x(N+1)=\xi$, we introduce a Lagrange multiplier $\lambda$ and define the following cost function:
\begin{eqnarray}
J_1&=&\sum_{k=0}^N[x'(k)Qx(k)+u'(k)Ru(k)]\nonumber\\
&&+x'(N+1)Hx(N+1)+2\lambda'x(N+1).\label{t40}
\end{eqnarray}
By applying the maximum principle \cite{MP}, there exists a unique solution to the optimal control of minimizing the cost function $J_1$ subject to (\ref{t1}) if and only if a kind of forward and backward equations has a unique solution where the forward equation is (\ref{t1}), the backward equation is governed by
\begin{eqnarray}
p(k-1)&=&A'(k)p(k)+Qx(k),\label{t7}
\end{eqnarray}
with terminal condition $p(N)=Hx(N+1)+\lambda$ and the optimal controller satisfies
\begin{eqnarray}
0&=&Ru(k)+B'(k)p(k).\label{t6}
\end{eqnarray}

Secondly, we solve the forward and backward equations (\ref{t1}), (\ref{t7}) and (\ref{t6}). The key is to verify that
\begin{eqnarray}
p(k)&=&P(k+1)x(k+1)+\eta(k).\label{t8}
\end{eqnarray}
where $P(k+1)$ is defined by (\ref{t3}) and $\eta(k)$ obeys
\begin{eqnarray}
\eta(k-1)&=A_c'(k)\eta(k),\label{t9}
\end{eqnarray}
with $\eta(N)=\lambda$.
In fact, at time $N$, (\ref{t8}) holds naturally since $P(N+1)=H$ and $\eta(N)=\lambda$.
Assume that (\ref{t8}) holds with $k=s$, that is, $p(s)=P(s+1)x(s+1)+\eta(s)$, we prove that (\ref{t8}) holds for $k=s-1$. By substituting (\ref{t8}) into (\ref{t6}), we have
\begin{eqnarray}
0&=&Ru(s)+B'(s)P(s+1)x(s+1)+B'(s)\eta(s)\nonumber\\
&=&[R+B'(s)P(s+1)B(s)]u(s)\nonumber\\
&&+B'(s)P(s+1)A(s)x(s)
+B'(s)\eta(s).\nonumber
\end{eqnarray}
In view of the fact that $R>0$ and $P(s+1)\geq 0$, it yields that $\Gamma(s)>0$. This implies that
\begin{eqnarray}
u(s)&=&K(s)x(s)-\Gamma^{-1}(s) B'(s)\eta(s).\label{t10}
\end{eqnarray}
By plugging (\ref{t10}) into (\ref{t1}), it is obtained that
\begin{eqnarray}
x(s+1)&=&A_c(s)x(s)-\bar{B}(s)\eta(s).\label{t11}
\end{eqnarray}
Combining with (\ref{t7}), we have
\begin{eqnarray}
p(s-1)&=&A'(s)P(s+1)A_c(s)x(s)\nonumber\\
&&-A'(s)P(s+1)\bar{B}(s)\eta(s)\nonumber\\
&&+A'(s)\eta(s)+Qx(s).\nonumber
\end{eqnarray}
(\ref{t8}) follows at time $k=s-1$ from (\ref{t3}) and (\ref{t9}).

Finally, we solve the Lagrange multiplier $\lambda$. From (\ref{t11}), it follows that
\begin{eqnarray}
x(N+1)&=&\Phi(0,N)x(0)-\sum_{k=0}^N\Phi(k+1,N)\bar{B}(k)\eta(k)\nonumber\\
&=&\Phi(0,N)x(0)-G(0)\lambda,\label{t12}
\end{eqnarray}
where $\eta(k)=\Phi'(k+1,N)\lambda$ has been used in the derivation of the last equality. Since $x(N+1)=\xi$, (\ref{t12})
is reduced to the equation (\ref{t4}). Thus, if (\ref{t4}) has a solution $\lambda^*$, then the optimal controller (\ref{t5}) follows from
(\ref{t10}). The proof is now completed. \hfill $\blacksquare$

From Lemma \ref{lem1},  one of the key to the optimal solution of the Problem is the solvability of the equation (\ref{t4}). To further derive the solvability condition of equation (\ref{t4}), we introduce the reachability of system (\ref{t1}) whose definition is given below.
\begin{definition}\label{d1}
The state $\xi$ is said to be reachable at time $N+1$ from initial value $x$ for system (\ref{t1}), if there exists a sequence of controllers $\{u(k), k=0, 1, \ldots, N\}$ such that $x(0)=x$ and $x(N+1)=\xi$.
\end{definition}
The criterion to judge the reachability is presented in the following lemma.
\begin{lemma}\label{lem3}
$\xi$ is reachable at time $N+1$ from the initial value $x$ for system (\ref{t1}) if and only if
$\xi-\Phi(0, N)x\in Range(G_1)$, that is, there exists $\zeta$, such that
\begin{eqnarray}
\xi-[\Pi_{k=N}^0 A(k)]x=G_1\zeta, \label{zeta}
\end{eqnarray}
where $G_1=\sum_{k=0}^{N}[\Pi_{i=N}^{k+1} A(i)]B(k)B'(k)[\Pi_{i=N}^{k+1} A(i)]'$, and $\Pi_{i=N}^{k} A(i)\triangleq A(N)\cdots\cdots A(k)$ with $\Pi_{i=N}^{N+1} A(i)\triangleq I$.

\end{lemma}
\emph{Proof.} ``Sufficiency" Assume that (\ref{zeta}) holds, we prove that $\xi$ is reachable from initial value $x$ for system (\ref{t1}).
In this case, we design the controller
\begin{eqnarray}
u(k)=B'(k)[\Pi_{i=N}^{k+1} A(i)]'\zeta,\nonumber
\end{eqnarray}
where $\zeta$ satisfies (\ref{zeta}).
Together with (\ref{t1}), it yields that
\begin{eqnarray}
&&x(N+1)\nonumber\\
&=&\Pi_{k=N}^0 A(k)x+\sum_{k=0}^{N}\Pi_{i=N}^{k+1} A(i)B(k)u(k)\nonumber\\
&=&\Pi_{k=N}^0 A(k)x+G_1\zeta\nonumber\\
&=&\xi.\nonumber
\end{eqnarray}
Thus, $\xi$ is reachable at time $N+1$ from the initial value $x$ for system (\ref{t1}).

``Necessity" Assume that $\xi$ is reachable at time $N+1$ from initial value $x$ for system (\ref{t1}), we prove that there exists a $\zeta$ such that (\ref{zeta}) holds. Otherwise, for any $\zeta$, it holds that $\xi-[\Pi_{k=N}^0 A(k)]x\neq G_1\zeta$.
Then there exists $\alpha$ such that $G_1\alpha=0$ and $\alpha'\{\xi-[\Pi_{k=N}^0 A(k)]x\}\neq0$.

On the one hand, from $G_1\alpha=0$, we have
\begin{eqnarray}
0=\alpha'G_1\alpha=\sum_{k=0}^{N}\|\alpha'\Pi_{i=N}^{k+1} A(i)B(k)\|^2, \nonumber
\end{eqnarray}
which implies that $\alpha'\Pi_{i=N}^{k+1} A(i)B(k)=0$ for $k=0, 1, \ldots, N$.

On the other hand, by denoting the admissible controller such that $x(0)=x, x(N+1)=\xi$ by $u(k)$, it is obtained from (\ref{t1}) that
\begin{eqnarray}
\xi-\Pi_{k=N}^0 A(k)x=\sum_{k=0}^{N}[\Pi_{i=N}^{k+1} A(i)]B(k)u(k).\nonumber
\end{eqnarray}
Combining with $\alpha'\{\xi-[\Pi_{k=N}^0 A(k)]x\}\neq0$, it follows that
\begin{eqnarray}
\sum_{k=0}^{N}\{\alpha'[\Pi_{i=N}^{k+1} A(i)]B(k)u(k)\}\neq 0,\nonumber
\end{eqnarray}
which is a contradiction to the fact that $\alpha'\Pi_{i=N}^{k+1} A(i)B(k)=0$.
Thus, there exists $\zeta$ such that (\ref{zeta}) holds. The proof is now completed. \hfill $\blacksquare$

We now show the solvability condition of equation (\ref{t4}).
\begin{lemma}
If $\xi$ is reachable at time $N+1$ from initial value $x$ for system (\ref{t1}), then the matrix equation (\ref{t4}) has a solution.
\end{lemma}
\emph{Proof.}
The proof is divided into two steps. The first step is to introduce a new system whose reachability is equivalent to that of (\ref{t1}).
The second step is to verify that the matrix equation (\ref{t4}) has a solution under the condition that the system (\ref{t1}) is reachable.

Firstly, under the reachability of $\xi$ at time $N+1$ from the initial value $x$ for system (\ref{t1}), there exists an admissible controller denoting by $u(k)$ such that $x(0)=x, x(N+1)=\xi$. We now introduce the following system:
\begin{equation}
\hat{x}(k+1)=A_c(k)\hat{x}(k)-\hat{B}(k)\hat{u}(k), \label{t1-hat}
\end{equation}
where the initial value is given by $\hat{x}(0)=x$ and
\begin{eqnarray}
\hat{B}(k)&=&B(k)[R+B'(k)P(k+1)B(k)]^{-\frac{1}{2}},\nonumber\\
\hat{u}(k)&=&-[R+B'(k)P(k+1)B(k)]^{\frac{1}{2}}u(k)\nonumber\\
&\ &-[R+B'(k)P(k+1)B(k)]^{-\frac{1}{2}}B'(k)P(k+1)\nonumber\\
&\ &\times A(k)x(k),\nonumber
\end{eqnarray}
while $x(k)$ is the solution of (\ref{t1}) with respect to the control $u(k)$.
Rewriting (\ref{t1-hat}) yields that
\begin{equation}
\hat{x}(k+1)=A(k)\hat{x}(k)+B(k)u(k)+B(k)K(k)[\hat{x}(k)-x(k)], \label{t1-hat-2}
\end{equation}
It is obvious that $x(k)$ satisfies (\ref{t1-hat}). Accordingly, under the controller $\hat{u}(k)$, the state $\hat{x}(k)$ in (\ref{t1-hat}) satisfies $\hat{x}(0)=x, \hat{x}(N+1)=\xi$. By using Lemma \ref{lem3}, there exists $\lambda_1$, such that
\begin{eqnarray}
\xi&=&\Phi(0,N)\hat{x}+G(0)\lambda_1, \label{x-N+1-hat}
\end{eqnarray}
that is, $\xi-\Phi(0,N)x=G(0)\lambda_1$. Thus, the matrix equation (\ref{t4}) has a solution $\lambda^*=-\lambda_1$.
The proof is now completed. \hfill $\blacksquare$

As a summary of Lemma \ref{lem1} and \ref{lem2}, we present the optimal solution to the Problem as follows.
\begin{theorem}\label{thm1}
If $\xi$ is reachable at time $N+1$ from initial value $x$ for system (\ref{t1}), then the equation (\ref{t4}) has a solution $\lambda^*$ and there exists a unique solution to the Problem which is given by (\ref{t5}).
\end{theorem}
\emph{Proof.} The proof follows directly by using Lemma \ref{lem1} and \ref{lem2}. The proof is now completed. \hfill $\blacksquare$

\section{Q-learning algorithm}

In this section, we aim to present the Q-learning algorithm for the Problem with unknown matrices $A(k)$ and $B(k)$.
From Theorem \ref{thm1}, the optimal controller (\ref{t5}) is in the feedback of the state and a specific parameter $\lambda^*$ which is the solution of the matrix equation (\ref{t4}). To this end, the design of the Q-learning algorithm is divided into two steps. The first step is to learn the feedback gain $K(k), K_1(k)$ and the solution $P(k)$ of the Riccati equation (\ref{t3}) as well the matrices $\Phi(k,N), G(k)$. This gives a data-driven representation of matrix gains $\Phi(0,N), G(0)$ in equation (\ref{t4}). Then, the second step is to solve the data-driven matrix equation (\ref{t4}) to derive the parameter $\lambda^*$.

\subsection{Q-function}

Firstly, to learn the feedback gain $K(k), K_1(k)$ and $P(k)$, we introduce the Q-function by using the dynamic programming. To this end, we first consider the optimal control problem of minimizing the cost function (\ref{t40}) with any parameter $\lambda$ and define the Q-function
for $s=0, \cdots, N$ by
\begin{eqnarray}
Q(s)&=&\sum_{k=s}^N[x'(k)Qx(k)+u'(k)Ru(k)]\nonumber\\
&&+x'(N+1)Hx(N+1)+2\lambda'x(N+1).\label{t13}
\end{eqnarray}
By applying the optimality principle, it follows that
\begin{eqnarray}
Q(s)&=&x'(k)Qx(k)+u'(k)Ru(k)+Q(s+1),\label{t14}
\end{eqnarray}
where $Q(N+1)=x'(N+1)Hx(N+1)+2\lambda'x(N+1)$. The representation of the Q-function can be given as below.
\begin{lemma}\label{lem2}
The Q-function defined by (\ref{t14}) is given by
\begin{eqnarray}
Q(s)=x'(s)P(s)x(s)+2x'(s)\Phi'(s,N)\lambda-\lambda'G(s)\lambda.\label{t15}
\end{eqnarray}
The optimal controller to minimize the cost function (\ref{t40}) with any parameter $\lambda$ is as
\begin{eqnarray}
u(s)&=&K(s)x(s)+K_1(s)\lambda.\label{t18}
\end{eqnarray}
\end{lemma}
\emph{Proof.} The proof is derived by using the induction technique. Firstly, from $s=N$ and using (\ref{t1}), it is obtained
\begin{eqnarray}
Q(N)&=&x'(N)Qx(N)+u'(N)Ru(N)\nonumber\\
&&+x'(N+1)Hx(N+1)+2\lambda'x(N+1)\nonumber\\
&=&x'(N)[A'(N)HA'(N)+Q]x(N)\nonumber\\
&&+u'(N)\Gamma(N)u(N)\nonumber\\
&&+2x'(N)A'(N)HB(N)u(N)+2\lambda'A(N)x(N)\nonumber\\
&&+2\lambda'B(N)u(N)\nonumber\\
&=&x'(N)P(N)x(N)\nonumber\\
&&+[u'(N)-K(N)x(N)-K_1(N)\lambda]'\nonumber\\
&&\times \Gamma(N)[u'(N)-K(N)x(N)-K_1(N)\lambda]\nonumber\\
&&+2x'(N)A_c'(N)\lambda-\lambda'\bar{B}(N)\lambda.\nonumber
\end{eqnarray}
Thus, the optimal controller is given by (\ref{t18}) and the corresponding Q-function is given by
\begin{eqnarray}
Q(N)
&=&x'(N)P(N)x(N)+2x'(N)\Phi'(N,N)\lambda\nonumber\\
&&-\lambda'G(N)\lambda.\nonumber
\end{eqnarray}
Assume that it holds that
\begin{eqnarray}
Q(s+1)&=&x'(s+1)P(s+1)x(s+1)\nonumber\\
&&+2x'(s+1)\Phi'(s+1, N)\lambda-\lambda'G(s+1)\lambda,\nonumber
\end{eqnarray}
we have that
\begin{eqnarray}
Q(s)&=&x'(s)Qx(s)+u'(s)Ru(s)\nonumber\\
&&+x'(s+1)P(s+1)x(s+1)\nonumber\\
&&+2x'(s+1)\Phi'(s+1, N)\lambda-\lambda'G(s+1)\lambda\nonumber\\
&=&x'(s)P(s)x(s)+[u(s)-K(s)x(s)\nonumber\\
&&-K_1(s)\lambda]'\Gamma(N)[u(s)\nonumber\\
&&-K(s)x(s)-K_1(s)\lambda]\nonumber\\
&&+2x'(s)\Phi'(s,N)\lambda-\lambda'G(s)\lambda.\nonumber
\end{eqnarray}
This implies that the optimal controller is given by (\ref{t18})
and the Q-function (\ref{t15}) follows.
The proof is now completed. \hfill $\blacksquare$

Based on Lemma \ref{lem2}, we now derive the data-driven representation of the Q-function.
\begin{eqnarray}
Q(s)
&=&x'(s)[Q+A'(s)P(s+1)A(s)]x(s)\nonumber\\
&&+u'(s)[R+B'(s)P(s+1)B(s)]u(s)\nonumber\\
&&+2x'(s)A'(s)P(s+1)B(s)u(s)\nonumber\\
&&+2x'(s)A'(s)\Phi'(s+1, N)\lambda\nonumber\\
&&+2u'(s)B'(s)\Phi'(s+1, N)\lambda-\lambda'G(s+1)\lambda\nonumber\\
&=&\left[
     \begin{array}{c}
       x(s) \\
       u(s) \\
       \lambda \\
     \end{array}
   \right]'\Lambda(s)\left[
     \begin{array}{c}
       x(s) \\
       u(s) \\
       \lambda \\
     \end{array}
   \right],\label{t16}
\end{eqnarray}
where
\begin{eqnarray}
\Lambda(s)&\triangleq& \left[
                         \begin{array}{ccc}
                           \Lambda_{11}(s) & \Lambda_{21}'(s) & \Lambda_{31}'(s) \\
                           \Lambda_{21}(s) & \Lambda_{22}(s) & \Lambda_{32}'(s) \\
                           \Lambda_{31}(s) & \Lambda_{32}(s) & \Lambda_{33}(s) \\
                         \end{array}
                       \right],\nonumber
                       \end{eqnarray}

\begin{eqnarray}
\Lambda_{11}(s)&=&Q+A'(s)P(s+1)A(s),\nonumber\\
\Lambda_{21}(s)&=&B'(s)P(s+1)A(s),\nonumber\\
\Lambda_{22}(s)&=&\Gamma(s),\nonumber\\
\Lambda_{31}(s)&=&\Phi(s+1, N)A(s),\nonumber\\
\Lambda_{32}(s)&=&\Phi(s+1, N)B(s),\nonumber\\
\Lambda_{33}(s)&=&-G(s+1).\nonumber
\end{eqnarray}

From (\ref{t16}), the feedback gains in (\ref{t18}) is reformulated as
\begin{eqnarray}
K(s)&=&-\Lambda_{22}^{-1}(s)\Lambda_{21}(s),\label{t17}\\
K_1(s)&=&-\Lambda_{22}^{-1}(s)\Lambda_{32}'(s).
\end{eqnarray}
The Riccati equation (\ref{t3}) is formulated by
\begin{eqnarray}
P(s)&=&\Lambda_{11}(s)-\Lambda_{21}'(s)\Lambda_{22}^{-1}(s)\Lambda_{21}(s).\label{t19}
\end{eqnarray}

\subsection{Q-learning algorithm}
We are now in the position to illustrate the derivation of the Q-learning algorithm. To this end, we take $l$ experiments by entering
the inputs $u^i(k)$, the state $x^i(k)$ where $i=1, 2, \ldots, l$,
and obtain the measurements $x^i(k+1)$ of the state for $k=0, 1, \ldots, N$. Moreover, we also enter parameter vector $\lambda^i(k), i=1, 2, \ldots, l$ to derive the data-driven Q-function (\ref{t15}) for $k=0, 1, \ldots, N$.

Firstly, with the acquired data information $(x^i(N), u^i(N), \lambda^i(N), x^i(N+1))$, we define
\begin{eqnarray}
\gamma^i(N)&=&[x^i(N)]'Qx^i(N)+[u^i(N)]'Ru^i(N)\nonumber\\
&&+[x^i(N+1)]'H[x^i(N+1)]\nonumber\\
&&+2[x^i(N+1)]'\lambda^i(N).\label{t35}
\end{eqnarray}
From (\ref{t16}), we have that
\begin{eqnarray}
\gamma^i(N)
&=&\left[
     \begin{array}{c}
       x^i(N) \\
       u^i(N) \\
       \lambda^i(N) \\
     \end{array}
   \right]'\Lambda(N)\left[
     \begin{array}{c}
       x^i(N) \\
       u^i(N) \\
       \lambda^i(N) \\
     \end{array}
   \right].
\end{eqnarray}
This gives the solvability of the matrix $\Lambda(N)$ as follows:
\begin{eqnarray}
\Lambda(N)=argmin\sum_{i=1}^l\|[z^i(N)]'\Lambda(N)z^i(N)-\gamma^i(N)\|,\label{t21}
\end{eqnarray}
where
\begin{eqnarray}
z^i(k)&=&\left[
     \begin{array}{c}
       x^i(N) \\
       u^i(N) \\
       \lambda^i(N) \\
     \end{array}
   \right]\nonumber\\
&\triangleq& \left[
                       \begin{array}{cccc}
                         z^i_1(N) &  z^i_2(N)& \cdots & z^i_{2n+m}(N) \\
                       \end{array}
                     \right].\nonumber
\end{eqnarray}
In fact, by denoting $\nu(N)$ the vectorization of the upper triangular part of
the matrix $\Lambda(N)$, the solvability of (\ref{t21}) can be further reduced to the following optimization problem:
\begin{eqnarray}
&&argmin \|\Upsilon(N)\nu(N)-\gamma(N)\|^2\nonumber\\
&=&[\Upsilon'(N)\Upsilon(N)]^{-1}\Upsilon'(N)\gamma(N),\label{t22}
\end{eqnarray}
where
\begin{eqnarray}
\Upsilon(k)&=&\left[
                \begin{array}{cccc}
                  \Upsilon^1(k) & \Upsilon^2(k)  & \cdots & \Upsilon^l(k)  \\
                \end{array}
              \right],\label{t39}\\
\gamma(k)&=&\left[
                \begin{array}{cccc}
                  \gamma^1(k) & \gamma^2(k)  & \cdots & \gamma^l(k)  \\
                \end{array}
              \right]',\nonumber
\end{eqnarray}
while
\begin{eqnarray}
\Upsilon^i(k)&=&\left[
                  \begin{array}{cccccccccc}
                    [z_{1}^i(k)]^2 & z_{1}^i(k)z_2^i(k) & \cdots & z_{1}^i(k)z_{2n+m}^i(k) &  & &  & && \\
                  \end{array}
                \right.\nonumber\\
&&\hspace{-16mm}\left.
                  \begin{array}{cccccccccc}
                    &&&&[z_{2}^i(k)]^2& z_{2}^i(k)z_3^i(k) & \cdots & z_{2}^i(k)z_{2n+m}^i(k) & \cdots &  \\
                  \end{array}
                \right.\nonumber\\
&&\hspace{-20mm}\left.
                  \begin{array}{cccccccccc}
                    &&&&&&&&& [z_{2n+m}^i(k)]^2 \\
                  \end{array}
                \right]'.\nonumber
\end{eqnarray}
Accordingly, the feedback gains at time $N$ are as
\begin{eqnarray}
K(N)&=&-\Lambda_{22}^{-1}(N)\Lambda_{21}(N),\label{t24}\\
K_1(N)&=&-\Lambda_{22}^{-1}(N)\Lambda_{32}'(N),\label{t25}
\end{eqnarray}
the solution $P(N)$ to the Riccati equation is given by
\begin{eqnarray}
P(N)&=&\Lambda_{11}(N)-\Lambda_{21}'(N)\Lambda_{22}^{-1}(N)\Lambda_{21}(N),\label{t23}
\end{eqnarray}
and the matrices $\Phi(N,N), G(N)$ are calculated by
\begin{eqnarray}
\Phi(N,N)&=&\Lambda_{31}(N)-\Lambda_{32}(N)\Lambda_{22}^{-1}(N)\Lambda_{21}(N),\label{t26}\\
G(N)&=&\Lambda_{32}(N)\Lambda_{22}^{-1}(N)\Lambda_{32}'(N).\label{t27}
\end{eqnarray}

Secondly, by repeatedly iterating backward from $N$, using the acquired data information $(x^i(k), u^i(k), \lambda^i(k), x^i(k+1))$ and the data-driven matrix representations in the last step, we define
\begin{eqnarray}
\gamma^i(k)&=&[x^i(k)]'Qx^i(k)+[u^i(k)]'Ru^i(k)\nonumber\\
&&+[x^i(k+1)]'P(k+1)x^i(k+1)\nonumber\\
&&+2[x^i(k+1)]'\Phi'(k+1,N)\lambda^i(k)\nonumber\\
&&-[\lambda^i(k)]'G(k+1)\lambda^i(k),\label{t37}
\end{eqnarray}
which can also be obtained from (\ref{t16}) that
\begin{eqnarray}
\gamma^i(k)
&=&\left[
     \begin{array}{c}
       x^i(k) \\
       u^i(k) \\
       \lambda^i(k) \\
     \end{array}
   \right]'\Lambda(k)\left[
     \begin{array}{c}
       x^i(k) \\
       u^i(k) \\
       \lambda^i(k) \\
     \end{array}
   \right].
\end{eqnarray}
Thus, $\Lambda(k)$ is solved as follows:
\begin{eqnarray}
\Lambda(k)=argmin\sum_{i=1}^l\|[z^i(k)]'\Lambda(k)z^i(k)-\gamma^i(k)\|,\label{t28}
\end{eqnarray}
Similarly to the derivation of (\ref{t22}), we solve (\ref{t28}) by the following optimization problem:
\begin{eqnarray}
&&argmin \|\Upsilon(k)\nu(k)-\gamma(k)\|^2\nonumber\\
&=&[\Upsilon'(k)\Upsilon(k)]^{-1}\Upsilon'(k)\gamma(k).\label{t29}
\end{eqnarray}
This gives that the feedback gains $K(k), K_1(k)$ is given by
\begin{eqnarray}
K(k)&=&-\Lambda_{22}^{-1}(k)\Lambda_{21}(k),\label{t30}\\
K_1(k)&=&-\Lambda_{22}^{-1}(k)\Lambda_{32}'(k),\label{t31}
\end{eqnarray}
the solution to the Riccati equation is calculated by
\begin{eqnarray}
P(k)&=&\Lambda_{11}(k)-\Lambda_{21}'(k)\Lambda_{22}^{-1}(k)\Lambda_{21}(k),\label{t32}
\end{eqnarray}
and the matrices $\Phi(k,N), G(k)$ are shown by
\begin{eqnarray}
\Phi(k,N)&=&\Lambda_{31}(k)-\Lambda_{32}(k)\Lambda_{22}^{-1}(k)\Lambda_{21}(k),\label{t33}\\
G(k)&=&G(k+1)+\Lambda_{32}(k)\Lambda_{22}^{-1}(k)\Lambda_{32}'(k).\label{t34}
\end{eqnarray}

Finally, we give the calculation of the parameter $\lambda^*$. From (\ref{t4}), we have
\begin{eqnarray}
\lambda^*=G^{\dag}(0)[\Phi(0,N)x-\xi].\nonumber
\end{eqnarray}
Combining with the data-driven representation in (\ref{t33}) and (\ref{t34}) at time $k=0$, it is further obtained that
\begin{eqnarray}
\lambda^*&=&[-\Lambda_{33}(0)+\Lambda_{32}(0)\Lambda_{22}^{-1}(0)\Lambda_{32}'(0)]^{\dag}\Big[[\Lambda_{31}(0)\nonumber\\
&&-\Lambda_{32}(0) \Lambda_{22}^{-1}(0)\Lambda_{21}(0)]x-\xi\Big].\label{t20}
\end{eqnarray}

The detailed algorithm is given in Algorithm \ref{alg2}.
\begin{algorithm}[h]
	\protect\protect\protect\protect\protect\protect\protect\caption{(Q-learning Algorithm)}
	\label{alg2} \begin{algorithmic}
		\STATE
\begin{enumerate}
  \item[{\color{blue}1:}] Implement $l$ times experiments on the system (\ref{t1}) where for the $i$-th experiment
   with $i=1, \ldots, l$, obtain the data samples $(x^i(k), u^i(k), \lambda^i(k), x^i(k+1))$ for $k=0, \ldots, N$.
  \item[{\color{blue}2:}] Calculate $\gamma^i(N)$ by (\ref{t35}) and solve the optimization problem (\ref{t22}). Calculate the feedback gains $K(N), K_1(N)$ by (\ref{t24})-(\ref{t25}), $P(N)$ by (\ref{t23}), the matrices $\Phi(N,N), G(N)$ by (\ref{t26})-(\ref{t27}).
  \item[{\color{blue}3:}] For $k=N-1$ to $k=0$, using repeated
steps to calculate $\gamma^i(k)$ by (\ref{t37}) and solve the optimization problem (\ref{t29}). Calculate the feedback gains $K(k), K_1(k)$ by (\ref{t30})-(\ref{t31}), $P(k)$ by (\ref{t32}), the matrices $\Phi(k,N), G(k)$ by (\ref{t33})-(\ref{t34}).
  \item[{\color{blue}4:}] Calculate the parameter $\lambda^*$ by (\ref{t20}).
  \item[{\color{blue}5:}] Derive the optimal controller $u(k)$ by (\ref{t5}).
\end{enumerate}
	\end{algorithmic}
\end{algorithm}

The convergence of Algorithm \ref{alg2} is given as follows.

\begin{theorem}\label{thm2}
Assume that $\xi$ is reachable at time $N+1$ from initial value $x$ for system (\ref{t1}).
Collect $l\geq \frac{1}{2}(2n+m)(2n+m+1)$ data samples $(x^i(k), u^i(k), \lambda^i(k)), i=1, 2, \ldots, l$ for $k=0, 1, \ldots, N$
which are Gaussian with arbitrary mean and arbitrarily positive-definite covariances, then the optimal solution of (\ref{t29}) exists and the controller (\ref{t5}) obtained in Algorithm 1 is the optimal solution to the Problem.
\end{theorem}
\emph{Proof.}  By collecting the data $x^i(k), u^i(k), \lambda^i(k), i=1, 2, \ldots, l$ which are Gaussian at time $k$,
we have that the experiments are linearly independent.
Considering that that the experiments times $l\geq \frac{1}{2}(2n+m)(2n+m+1)$,
it is known that the formed data matrix $\Upsilon(k)$ in (\ref{t39}) has full rank \cite{Gaussian},
which implies that (\ref{t29}) holds. Together with Theorem \ref{thm1}, the derived controller (\ref{t5}) in Algorithm 1 is optimal for the Problem. The proof is now completed. \hfill $\blacksquare$

\section{Numerical Examples}

In this section, we present an example to verify the effectiveness of the proposed algorithm.
Consider system (\ref{t1}) with parameters given by
\begin{eqnarray}
&&A(0)=\left[
    \begin{array}{cc}
      1 & 2 \\
      -1 & 4 \\
    \end{array}
  \right], A(1)=\left[
    \begin{array}{cc}
      5 & 3 \\
      -2 & 1 \\
    \end{array}
  \right], \nonumber\\
&&A(2)=\left[
    \begin{array}{cc}
      -4 & 1 \\
      2 & 5 \\
    \end{array}
  \right], B(0)=\left[
               \begin{array}{c}
                 1 \\
                 -1 \\
               \end{array}
             \right], B(1)=\left[
               \begin{array}{c}
                 2 \\
                 1 \\
               \end{array}
             \right], \nonumber\\
&&B(2)=\left[
               \begin{array}{c}
                  4 \\
                  2 \\
               \end{array}
             \right], Q=I, R=1, N=2,\nonumber\\
&&x=\left[
               \begin{array}{c}
                  1 \\
                  2 \\
               \end{array}
             \right], \xi=\left[
               \begin{array}{c}
                  6 \\
                  7 \\
               \end{array}
             \right].\nonumber
\end{eqnarray}
By implementing $l=30$ times experiments, we have the solution of the optimization problem (\ref{t29}) as
\begin{eqnarray}
\nu(2)&=&\left[
           \begin{array}{c}
             21 \\
             6 \\
             -12 \\
             -4 \\
             2 \\
             27 \\
             14 \\
             1 \\
             5 \\
             21 \\
             4 \\
             2 \\
             0 \\
             0 \\
             0 \\
           \end{array}
         \right], \nu(1)=\left[
           \begin{array}{c}
             145.2381 \\
             162.8095\\
             120.0952 \\
             -5.2381 \\
             8.3810 \\
             229.9524 \\
             172.5238 \\
             -6.8095 \\
             13.0952 \\
             131.2381 \\
             -5.0952 \\
             9.9524  \\
             -0.7619 \\
             -0.3810 \\
             -0.1905 \\
           \end{array}
         \right], \nonumber\\
         &&\nu(0)=\left[
           \begin{array}{c}
            29.6266 \\
            67.9274\\
            28.6266 \\
            -0.4641 \\
            -0.7384 \\
            271.7808 \\
            67.9274 \\
            -1.5965 \\
            -1.4049 \\
            29.6266 \\
            -0.4641 \\
            -0.7384  \\
            -0.9597 \\
            0.0054 \\
            -0.9452 \\
           \end{array}
         \right].\nonumber
\end{eqnarray}
This accordingly gives the solution of the Riccati equation (\ref{t32}) as
\begin{eqnarray}
&&P(3)=I, P(2)=\left[
               \begin{array}{cc}
                 14.1429 & 14 \\
                 14 & 17.6667 \\
               \end{array}
             \right], \nonumber\\
&&P(1)=\left[
               \begin{array}{cc}
                 35.3396 & 4.9340 \\
                 4.9340 & 3.1549 \\
               \end{array}
             \right], \nonumber\\
&&P(0)=\left[
               \begin{array}{cc}
                 1.9662 & 2.2928 \\
                 2.2928 & 116.0380 \\
               \end{array}
             \right],\nonumber
\end{eqnarray}
and the feedback gains (\ref{t30})-(\ref{t31}) are given by
\begin{eqnarray}
K(0)&=&\left[
       \begin{array}{cc}
         -0.9662 & -2.2928 \\
       \end{array}
     \right],\nonumber\\
 K(1)&=&\left[
       \begin{array}{cc}
         -0.9151 & -1.3146 \\
       \end{array}
     \right],\nonumber\\
K(2)&=&\left[
       \begin{array}{cc}
         0.5714  & -0.6667 \\
       \end{array}
     \right], \nonumber\\
K_1(0)&=&\left[
       \begin{array}{cc}
         0.0157  & 0.0249 \\
       \end{array}
     \right],\nonumber\\
K_1(1)&=&\left[
       \begin{array}{cc}
         0.0388  & -0.0758 \\
       \end{array}
     \right], \nonumber\\
K_1(2)&=&\left[
       \begin{array}{cc}
         -0.1905  & -0.0952 \\
       \end{array}
     \right],\nonumber
\end{eqnarray}
and the parameter (\ref{t20}) is given by
\begin{eqnarray}
\lambda^*=\left[
            \begin{array}{c}
              -7.2802 \\
              -6.6461 \\
            \end{array}
          \right].\nonumber
\end{eqnarray}
It is easily verified that the closed-loop system (\ref{t1}) under the optimal controller (\ref{t5}) achieves the specific terminal state $\xi=\left[
\begin{array}{cc}
6 & 7 \\
\end{array}
\right]'$.

\section{Conclusions}

In this paper, we considered the linear quadratic optimal control of discrete-time time-varying system with terminal state constraint. A Q-learning algorithm was proposed for the optimal controller with unknown time-varying system matrices and input matrices, consisting of the learning of the solution to the Riccati equation and the calculation of the parameter solution to a data-driven matrix equation. Numerical example verified the effectiveness of the proposed algorithm.


\end{document}